\documentclass[]{interact}

\usepackage{epstopdf}
\usepackage[caption=false]{subfig}

\usepackage[numbers,sort&compress]{natbib}
\usepackage{tikz}

\theoremstyle{plain}

\theoremstyle{definition}

\theoremstyle{remark}

\newcommand\E{{\mathbb E}}
\newcommand\Poi{{\rm Poi}}
\newcommand\Prob{{\mathbb P}}
\newcommand\Exp{{\rm Exp}}
\newcommand{\almostsure}{\, \overset{a.s.}{\longrightarrow} \,}
\newcommand{\Polya}{P\'{o}lya}
\newcommand{\inL}{\, \overset{L_1}{\longrightarrow} \,}
\begin{document}

\articletype{ARTICLE TEMPLATE}

\title{The degree Gini index of several classes of random trees and their poissonized counterparts---an evidence for a duality theory}

\author{
\name{Carly Domicolo\textsuperscript{a}\thanks{CONTACT P.~Zhang. Email: panpan.zhang@pennmedicine.upenn.edu}, Panpan Zhang\textsuperscript{b} and Hosam Mahmoud\textsuperscript{a}}
\affil{\textsuperscript{a}Department of Statistics,	The George Washington University, Washington, D.C. 20052, U.S.A.; \textsuperscript{b}Department of Biostatistics, Epidemiology and Informatics, Perelman School of Medicine,
	University of Pennsylvania, Philadelphia, PA 19104, U.S.A.}
}

\maketitle

\begin{abstract}
There is an unproven duality theory hypothesizing that random discrete trees and their poissonized embeddings in continuous time share fundamental properties. We give additional evidence in favor 
of this theory by showing that several classes of random trees growing in discrete time and their poissonized counterparts
have the same limiting degree Gini index. The classes that we consider include binary search trees, binary pyramids 
and random caterpillars.
\end{abstract}
\section{Introduction}
There is an unproven duality theory hypothesizing that random discrete trees and 
certain embeddings in continuous time share fundamental properties. The embedding is done by changing the interarrival times between nodes from equispaced discrete units to more general renewal intervals. The most common method of embedding uses exponential random variables as the interarrival time. 
Under this choice of interarrival time, a count of the points of arrival constitutes a Poisson process. Hence, this kind of embedding is commonly called {\em poissonization}~\cite{Aldous, Sparks}. 
The advantage of poissonization is that the underlying exponential interarrival times have memoryless properties, giving rise to tractable features in the poissonized random structures. 

Let us call the discrete time $n$ and the continuous time $t$. The duality theory, now almost a folklore in the probability community, is a claim that what happens at discrete time $n$ has an analogue in the poissonized counterpart at continuous time $g(t)$, for some real function $g$. For instance, when the number of insertion positions in a discrete-time tree grows linearly in $n$, a poissonized tree has similar features at time $g(t) = e^t$; when the number of insertion positions in a discrete-time remains fixed as~$n$ grows, a poissonized tree has similar features at time $g(t) = t$.

Our intent in this paper is to substantiate the duality theory through a recently-developed topological index---the {\em Gini index} based on node degrees~\cite{Domicolo}. Indeed, in several classes of random trees growing in discrete time, their asymptotic (degree) Gini index is the same as the Gini index of the poissonized counterparts,
 evidence lending additional credence to the validity of the duality theory. 
\section{Gini index of a tree}
A topological index of a graph is a descriptor that quantifies its structure or 
some features
it possesses. Such quantification allows for practical comparison of graphs from different classes according to certain criteria. There are many indices that can be constructed for trees.
Each index captures a distinct feature of the graphs, 
 such as sparseness, thinness, regularity, centrality, etc. Examples of indices that have been introduced for trees include Zagreb index~\cite{Feng2011}, Randi\'{c} index~\cite{Feng2010}, Weiner index~\cite{Neininger}, a (distance-based) Gini index~\cite{Balaji}, two kinds of (degree-based) Gini index respectively in~\cite{Domicolo} and in~\cite{Zhang}. 

In the rest of this manuscript, we place our focus on the degree-based Gini indices
in the last two sources. In both, 
the degree Gini index is considered a measure of the degree profile of a tree, but with different applications. The degree-based Gini index in~\cite{Domicolo} is exploited to assess graph regularity, while the authors of~\cite{Zhang} use it to evaluate the balance of the degree distribution in random trees evolving in different manners. We will elaborate on this point in the sequel.
\section{Overview of the standard Gini index} 
\label{Sec:Gini}
The {\em Gini index} (or {\em Gini coefficient})~\cite{Gini} is a statistical measure of inequality, commonly used in 
economics~\cite{Gastwirth} for independent and identically distributed (i.i.d.) observations on income and wealth. A modified version of the Gini index suited for graphs 
is introduced in~\cite{Domicolo}. The computations in that index are based on the degrees of nodes in the graph. In random graphs, the node degrees are not i.i.d.\ in general, so the definition of Gini index in~\cite{Domicolo} imposes several tweaks on the standard Gini index.

Independently, the authors of~\cite{Zhang} consider a different version of degree-based Gini index, specifically designed for a
class of random trees called {\em caterpillars}, known for applications in 
chemistry and physics; see~\cite{ElBasil} for instance. 
Mostly, we shall use the definition in~\cite{Domicolo} for computations in various classes of trees. 
For caterpillars, we compare the discrete-time and poissonized versions via the definition in~\cite{Zhang}, too. 
There are two classes of caterpillars, respectively growing in a uniform and a nonuniform manner, considered in this manuscript.

To cast a Gini index 
as a topological measure on graphs, we first establish some useful notation. Suppose we have a graph $H = (V, E)$, where~$V$ and $E$ are the sets of {\em vertices} (nodes) 
and {\em edges}, respectively. 
Let $|V|$ be the {\em cardinality} of set $V$, and $|V|$ is also called the {\em order} of graph. For a graph $H$ of order $|V| = n$, we arbitrarily label the distinct vertices in $V$ with a set of distinct integers, $[n] := \{1, 2, \ldots, n\}$. Let $D_i(H)$ denote the degree of the node labeled with $i \in [n]$ in $H$, and let $D^*(H)$ be the degree of a node uniformly chosen from $V$, i.e., we set
$$\Prob\bigl(D^*(H) = D_i(H)\bigr) = 1/n,$$
for $i \in [n]$.
Therefore, the average node degree is given by
$$ \E \bigl[ D^*(H) \bigr] =\frac{1}{n} \sum_{i = 1}^n D_i(H).$$
The degree-based Gini index of $H$ in~\cite{Domicolo} is defined analogously to the standard Gini index, with the utilization of graph parameters; 
that is,  
\begin{equation}
G(H) = \frac {\sum_{1 \le i \le j \le n} \bigl|D_j(H)- D_i(H)\bigr|}{n^2 \,  \E \bigl[D^* (H) \bigr]}.
\label{Eq:singlegraph}
\end{equation}
For simplicity, we call this measure the {\em degree Gini index} of graph $H$. This is a topological index of a graph.

The definition (c.f.\ Equation~(\ref{Eq:singlegraph})) can be extended to a class of random graphs. Suppose that we have a class of random graphs, denoted 
by $\mathcal{H}$. For each graph $H \in \mathcal{H}$, one can define a {\em relative degree-based Gini index} for random $H$ by mimicking the topological degree Gini index in Equation~(\ref{Eq:singlegraph}). In contrast to the topological degree Gini index, the relative degree-based Gini index of $H$ 
measures the sum of the absolute degree differences in $H$ against global parameters taken from the entire class $\mathcal{H}$. The order of a graph $H\in \mathcal{H}$ is $|V(H)|$. Let $D^*_\mathcal{\mathcal{H}}$ be the degree of a randomly selected node in a randomly selected graph from $\mathcal{H}$. The relative degree-based Gini index of $H \in \mathcal{H}$ is given by
\begin{equation}
G_{\cal H}(H) = \frac {\sum_{1 \le i \le j \le |V(H)|} \bigl|D_j(H)- D_i(H)\bigr|}  {\E^2 |V(H)| \,  \E [D^*_{\cal H}]}.
\label{Eq:class}
\end{equation}
Taking the average over all $H \in \mathcal{H}$, 
we 
obtain the degree-based Gini index of the class~$\cal H$, which is namely
$$ G_{\cal H}^* = \E \bigl[G_{\cal H}(H)\bigr].$$ 

Note that if the class $\cal H$ has only one graph in it, the degree-based Gini index for $\mathcal{H}$ is then reduced to the topological degree Gini index of that one graph. A few illustrative examples are presented in~\cite{Domicolo}. We refer the interested readers to that article for more details. In the rest of the manuscript, we investigate the degree-based Gini index for several random classes, in which the random graphs evolve both in discrete time and in poissonized embeddings, in support of the duality theory. We shall see that in the chosen classes, the limiting degree Gini indices 
in the discrete-time trees and their poissonized counterparts coincide.

\section{Binary trees}
A {\em binary tree} is a structure of a finite number of nodes and edges. One special node is recognized as the {\em root}. The tree is either empty or has nodes, and each node can have at most two children. 
Thus, a binary tree is either empty or has a root and two subtrees (at the root) that are themselves binary trees. 

The binary tree is a well-studied data structure for the primary purpose of data storage and sorting~\cite[Section 6.1]{Sedgewick}. Many advanced
methods have been developed to analyze the properties of binary trees, such as generating functions~\cite{Sedgewick}, analytic combinatorics~\cite{Flajolet} and graph percolation~\cite{Mahmoudbook}. 

It expedites the analysis to consider an {\em extended binary tree} which can be thought of as a transformation of binary tree into a {\em full} binary tree, with uniform outdegrees. An extended binary tree is obtained from a given binary tree by supplying each node with a sufficient number of nodes of a special type called {\em external} to render the number of children of each node in the original tree exactly equal to 2. 

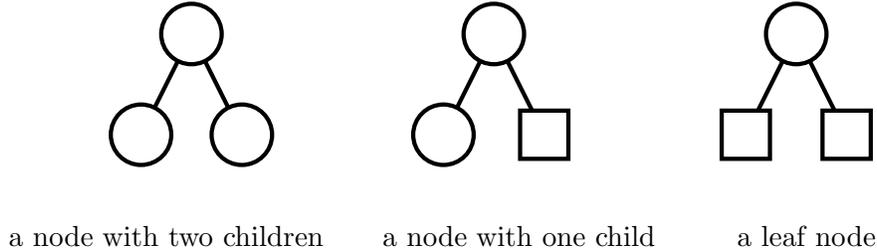
\begin{figure}[thb]
	\begin{center}
		\begin{tikzpicture}[scale=0.67]
		\coordinate (A) at (-3,4);
		\coordinate (B) at (3,4);
		\coordinate (C) at (9,4);
		
		\coordinate (D) at (-4,2);
		\coordinate (E) at (-2,2);
		\draw [ultra thick] (A)--(D);
		\draw [ultra thick] (A)--(E);
		\draw [ultra thick,fill=white] (-3,4) circle [radius=0.6];
		\draw [ultra thick,fill=white] (-4,2) circle [radius=0.6];
		\draw [ultra thick,fill=white] (-2,2) circle [radius=0.6];
		
		\coordinate (F) at (2,2);
		\coordinate (G) at (4,2);
		\draw [ultra thick] (B)--(F);
		\draw [ultra thick] (B)--(G);
		\draw [ultra thick,fill=white] (3,4) circle [radius=0.6];
		\draw [ultra thick,fill=white] (2,2) circle [radius=0.6];
		\node at (4,2) [rectangle, draw, ultra thick, scale = 2.5, fill=white](){};
		
		\coordinate (H) at (8,2);
		\coordinate (I) at (10,2);
		\draw [ultra thick] (C)--(H);
		\draw [ultra thick] (C)--(I);
		\draw [ultra thick,fill=white] (C) circle [radius=0.6];
		\node at (8,2) [rectangle, draw, ultra thick, scale = 2.5, fill=white](){};
		\node at (10,2) [rectangle, draw, ultra thick, scale = 2.5, fill=white](){};
		
		\node[draw=white] at (-3.5, 0) {\mbox{a node with two children}};
		\node[draw=white] at (3.5 , 0) {\mbox{a node with one child}};
		\node[draw=white] at (9.2, 0) {\mbox{a leaf node}};
		\end{tikzpicture}
	\end{center}
	\caption{The representation of three types of nodes in extended binary trees; round nodes are internal and square nodes are external.}
	\label{Fig:Extendedbinary}
\end{figure}

A model of randomness takes the external nodes as equally likely candidates for the next insertion. The model results in a random binary search tree. There are three types on internal nodes in the original tree---leaf nodes carrying no child, nonleaf nodes carrying one child, and nonleaf nodes carrying  two children. 
Those internal nodes having two children are {\em unsaturated}. The extension of these three types of nodes in binary trees is illustrated in Figure~\ref{Fig:Extendedbinary}.

\subsection{Gini index of discrete-time binary search trees}

Let $B(n)$ be a binary tree on $n$ nodes, of which there are $n_1$ nodes of degree $1$ (the leaves), $n_2$ nodes of degree $2$, and $n_3$ nodes of degree $3$, such that $n_1 + n_2 + n_3 = n$. The topological degree Gini index of $B(n)$, according to Equation~(\ref{Eq:singlegraph}), is
\begin{align}
G\bigl(B(n)\bigr) 
= \frac {n_1 n_2 + n_2n_3 + 2 n_1 n_3} {n^2 \times (n_1 + 2n_2 + 3n_3)/n} = \frac {n_1 n_2 + n_2n_3 + 2 n_1 n_3} {n \times (n_1 + 2n_2 + 3n_3)}.
\label{Eq:binGini}
\end{align}

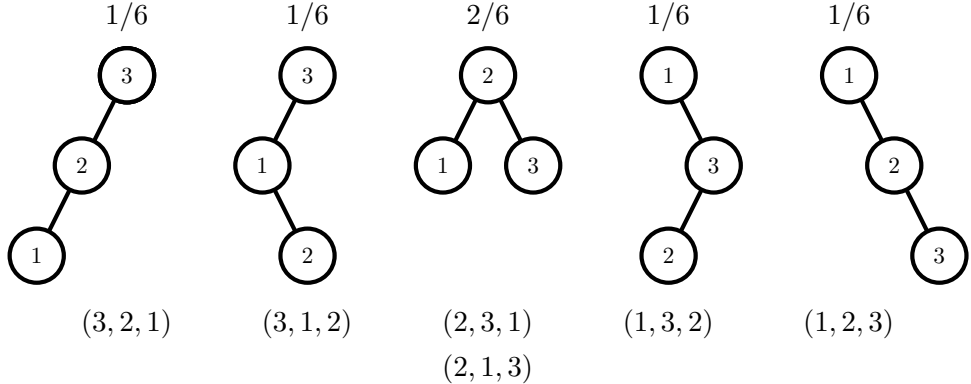
\begin{figure}[thb]
	\begin{center}
		\begin{tikzpicture}[scale=0.6]
		\coordinate (A) at (-5,4);
		\coordinate (B) at (-1,4);
		\coordinate (C) at (3,4);
		\coordinate (D) at (7,4);
		\coordinate (E) at (11,4);
		
		\coordinate (F) at (-6,2);
		\coordinate (G) at (-7,0);
		\draw [ultra thick] (A)--(F)--(G);
		\draw [ultra thick,fill=white] (-6,2) circle [radius=0.6];
		\draw [ultra thick,fill=white] (-7,0) circle [radius=0.6];

		\coordinate (H) at (-2,2);
		\coordinate (I) at (-1,0);
		\draw [ultra thick] (B)--(H)--(I);
		\draw [ultra thick,fill=white] (-2,2) circle [radius=0.6];
		\draw [ultra thick,fill=white] (-1,0) circle [radius=0.6];

		\coordinate (J) at (2,2);
		\coordinate (K) at (4,2);
		\draw [ultra thick] (K)--(C)--(J);
		\draw [ultra thick,fill=white] (2,2) circle [radius=0.6];
		\draw [ultra thick,fill=white] (4,2) circle [radius=0.6];
		
		\coordinate (L) at (8,2);
		\coordinate (M) at (7,0);
		\draw [ultra thick] (M)--(L)--(D);
		\draw [ultra thick,fill=white] (8,2) circle [radius=0.6];
		\draw [ultra thick,fill=white] (7,0) circle [radius=0.6];
		
		\coordinate (N) at (12,2);
		\coordinate (O) at (13,0);
		\draw [ultra thick] (E)--(N)--(O);
		\draw [ultra thick,fill=white] (12,2) circle [radius=0.6];
		\draw [ultra thick,fill=white] (13,0) circle [radius=0.6];

		\draw [ultra thick,fill=white] (-5,4) circle [radius=0.6];
		\draw [ultra thick,fill=white] (-1,4) circle [radius=0.6];
		\draw [ultra thick,fill=white] (3,4) circle [radius=0.6];
		\draw [ultra thick,fill=white] (7,4) circle [radius=0.6];
		\draw [ultra thick,fill=white] (11,4) circle [radius=0.6];

		\draw [ultra thick,fill=white] (-5,4) circle [radius=0.6];
		\draw [ultra thick,fill=white] (-5,4) circle [radius=0.6];
		
		\node[draw=white, scale = 0.8] at (-5, 4) {$3$};
		\node[draw=white, scale = 0.8] at (-6, 2) {$2$};
		\node[draw=white, scale = 0.8] at (-7, 0) {$1$};
		
		\node[draw=white, scale = 0.8] at (-1, 4) {$3$};
		\node[draw=white, scale = 0.8] at (-2, 2) {$1$};
		\node[draw=white, scale = 0.8] at (-1, 0) {$2$};
		
		\node[draw=white, scale = 0.8] at (3, 4) {$2$};
		\node[draw=white, scale = 0.8] at (2, 2) {$1$};
		\node[draw=white, scale = 0.8] at (4, 2) {$3$};

		\node[draw=white, scale = 0.8] at (7, 4) {$1$};
		\node[draw=white, scale = 0.8] at (8, 2) {$3$};
		\node[draw=white, scale = 0.8] at (7, 0) {$2$};

		\node[draw=white, scale = 0.8] at (11, 4) {$1$};
		\node[draw=white, scale = 0.8] at (12, 2) {$2$};
		\node[draw=white, scale = 0.8] at (13, 0) {$3$};
		
		\node[draw=white] at (-5, 5.3) {$1/6$};
		\node[draw=white] at (-1, 5.3) {$1/6$};
		\node[draw=white] at (3, 5.3) {$2/6$};
		\node[draw=white] at (7, 5.3) {$1/6$};
		\node[draw=white] at (11, 5.3) {$1/6$};
		
		\node[draw=white] at (-5, -1.5) {$(3,2,1)$};
		\node[draw=white] at (-1, -1.5) {$(3,1,2)$};
		\node[draw=white] at (3, -2.5) {$(2,1,3)$};
		\node[draw=white] at (3, -1.5) {$(2,3,1)$};
		\node[draw=white] at (7, -1.5) {$(1,3,2)$};
		\node[draw=white] at (11, -1.5) {$(1,2,3)$};
		
		\end{tikzpicture}
	\end{center}
	\caption{Binary search trees of order 3 arising from permutations of $\{1,2,3\}$. Below each tree are the permutations associated with it, and above each is its probability.}
	\label{Fig:BST}
\end{figure}

A popular probability measure on binary trees is induced by random permutations, giving rise to {\em binary search trees} (BST's)~\cite[Section 6.6]{Sedgewick}. As mentioned, binary trees are utilized as data storage devices. A model of relevance to the use of binary trees as data storage devices is their growth from a random permutation. Suppose that $K_1, \ldots, K_n$ are keys sampled from a continuous distribution. Their ranks then almost surely are a random 
permutation of $\{1, \ldots, n\}$. A binary tree is created for the storage
of these keys. The key $K_1$ is retained in a root node, with empty left and right subtrees. When a subsequent key appears, it is guided to the left or right subtree according as whether it is less than $K_1$ or is at least as large. The ranks are the only facet of the data that drives the shape of the resulting tree. The ranks are almost surely a random permutation of $\{1, \ldots, n\}$. The process then can be assimilated by the insertion of the~$n$ distinct elements of a random permutation of $\{1, \ldots, n\}$. The resulting tree is a binary search tree (BST). In Figure~\ref{Fig:BST}, we depict the five different shapes of BST's arising from random permutations of $\{1, 2, 3\}$. The probability of each shape is shown at the top. Apparently, the growth of BST's from permutations gives rise to a nonuniform distribution on the possible binary tree shapes.  

We call the class of binary search trees on $n$ nodes by the name
${\cal B}(n)$. Let $N_1$, $N_2$ and $N_3$ be the number of nodes of degree $1$, $2$, and $3$, respectively. It is shown in~\cite{Mahmoud1986, Devroye} that $N_1/n \to 1/3$, $N_2/n \to 1/3$ and $N_3/n \to 1/3$, where 
all convergence relations take place in the $L_1$ space as well as in probability. We refer the readers to~\cite[Section 8.1]{Mahmoud2008} for text-style elaboration of these results. As in~\cite{Domicolo}, we obtain the degree-based Gini index for the class of binary search trees:
$$  G^*_ {{\cal B}(n) }  \to \frac{2}{9},$$
as $n \to \infty$.
\subsection{Gini index of poissonized binary search trees}
We focus next on the degree-based Gini index of poissonized BST's. We denote an exponential random varaible with mean $\lambda > 0$ by $\Exp(\lambda)$. In this model of growth, we have an independent
$\Exp(1)$ random variable associated with each insertion position (all 
random exponential random variables being independent) in an extended BST, as if a percolating fluid is coming down the edges connecting unsaturated internal nodes to external nodes, at an independent speed, each finishing its percolation in an independent $\Exp(1)$ time. Almost surely, one of the racing droplets will win the competition and be the first to touch the external node at the end of its edge. At this point, the wetted external node is converted into a leaf node, out of which two edges lead to two new external nodes and each new edge is endowed with an independent $\Exp(1)$ random variable. By the memoryless property of exponential random variables, the incomplete percolation process on all the other edges 
is reset. 
In other words, right after an insertion, the remaining time to finish the percolation on each edge remains $\Exp(1)$, and that includes the two new edges, too. This analogy is borrowed from~\cite{Pittel}. At time $t$, we call the class of trees 
generated ${\cal B}_{\mathsf P}(t)$.

This model of growth corresponds to a two-color \Polya\ urn~\cite{Mahmoud2008} evolving in real time. The insertion positions (external nodes) under leaf nodes are the white balls of the urn, and the remaining insertion positions (under internal nodes already having one child) are the blue balls of the urn. 

\begin{figure}[tbh]
	\begin{center}
		\begin{minipage}{0.48\textwidth}
			\begin{tikzpicture}[scale=0.67]
			\coordinate (A) at (-3,4);
			\coordinate (B) at (3,4);
			\node[draw=white] at (-4, 2.8) {$\star$};
			\coordinate (D) at (-4,2);
			\coordinate (E) at (-2,2);
			\draw [ultra thick] (A)--(D);
			\draw [ultra thick] (A)--(E);
			\draw [ultra thick,fill=white] (-3,4) circle [radius=0.6];
			\node at (-4,2) [rectangle, draw, ultra thick, scale = 2.5, fill=white](){};
			\node at (-2,2)[rectangle, draw, ultra thick, scale = 2.5, fill=white](){};
			
			\draw [->,line width=1pt] (-0.5,2) -- (0.5,2); 
			
			\coordinate (F) at (2,2);
			\coordinate (G) at (4,2);
			\coordinate (H) at (1,0);
			\coordinate (I) at (3,0);
			\draw [ultra thick] (B)--(F);
			\draw [ultra thick] (B)--(G);
			\draw [ultra thick] (F)--(H);
			\draw [ultra thick] (F)--(I);
			\draw [ultra thick,fill=white] (3,4) circle [radius=0.6];
			\draw [ultra thick,fill=white] (2,2) circle [radius=0.6];
			\node at (4,2) [rectangle, draw, ultra thick, scale = 2.5, fill=blue](){};
			\node at (1,0) [rectangle, draw, ultra thick, scale = 2.5, fill=white](){};
			\node at (3,0) [rectangle, draw, ultra thick, scale = 2.5, fill=white](){};
			\end{tikzpicture}
		\end{minipage}
	\hfill
		\begin{minipage}{0.48\textwidth}
			\begin{tikzpicture}[scale=0.67]
					\coordinate (A) at (-3,4);
					\coordinate (B) at (3,4);
					\node[draw=white] at (-2, 2.8) {$\star$};
					\coordinate (D) at (-4,2);
					\coordinate (E) at (-2,2);
					\draw [ultra thick] (A)--(D);
					\draw [ultra thick] (A)--(E);
					\draw [ultra thick,fill=white] (-3,4) circle [radius=0.6];
					\draw [ultra thick,fill=white] (-4,2) circle [radius=0.6];
					\node at (-2,2)[rectangle, draw, ultra thick, scale = 2.5, fill=blue](){};
					
					\draw [->,line width=1pt] (-0.5,2) -- (0.5,2); 
					
					\coordinate (F) at (2,2);
					\coordinate (G) at (4,2);
					\coordinate (H) at (3,0);
					\coordinate (I) at (5,0);
					\draw [ultra thick] (B)--(F);
					\draw [ultra thick] (B)--(G);
					\draw [ultra thick] (G)--(H);
					\draw [ultra thick] (G)--(I);
					\draw [ultra thick,fill=white] (3,4) circle [radius=0.6];
					\draw [ultra thick,fill=white] (2,2) circle [radius=0.6];
					\draw [ultra thick,fill=white] (4,2) circle [radius=0.6];
					\node at (3,0) [rectangle, draw, ultra thick, scale = 2.5, fill=white](){};
					\node at (5,0) [rectangle, draw, ultra thick, scale = 2.5, fill=white](){};
					\end{tikzpicture}
		\end{minipage}
	\caption{The ball evolution of two-color P\'{o}lya urns in an extended BST.}
	\label{Fig:polyaurn}
	\end{center}
\end{figure}
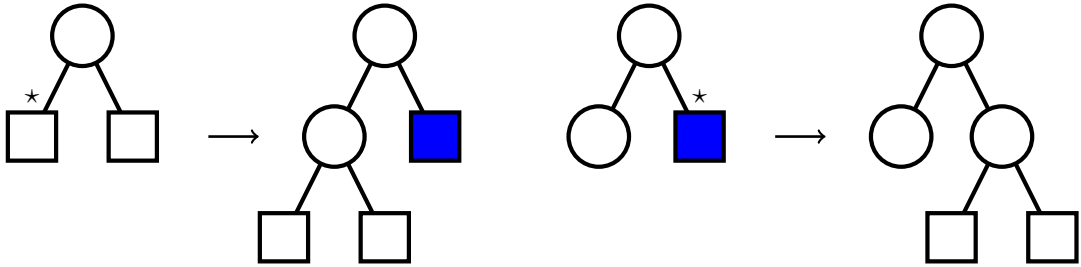

By the memoryless property of exponential random variables, at any point in time, each edge involving external nodes carries an independent $\Exp(1)$. With the exponential variables representing 
the growth along the edges leading to these insertion positions being i.i.d., any of them can achieve the minimum interarrival time. That is to say, all the balls in the urn are equally likely to be picked. When a white ball is picked, the corresponding insertion position is converted into a leaf node, carrying two new leaves
(corresponding to two white balls in the urn). A blue ball appears 
to represent a sibling. So, there is no change in the number of white balls, and the net gain of blue balls is one. Alternatively, if a blue ball is picked, the corresponding insertion position is converted into a leaf node, with two leaves underneath it  (two white balls in the urn). The total gain is two white balls, and the total loss is one blue ball. 
See Figure~\ref{Fig:polyaurn} for a depiction of ball evolution. The external nodes chosen for insertion are starred.
 
The ball addition matrix, also called {\em replacement matrix}, is
$${\bf A} = \begin{pmatrix} 0  &1\\
2   &-1
\end{pmatrix};$$
the rows correspond to the colors of the balls picked, 
with white at the top and blue at the bottom, and the 
columns correspond to the balls added, with white on the left and blue on the right. The transpose of $\bf A$ will be indicated by ${\bf A}^\mathsf{T}$.
 
Let $W(t)$ be the number of white balls in the urn at time $t$, and $B(t)$ be the number of blue balls in the urn at time $t$.
According to a functional urn theory~\cite{Janson, Mahmoud2008}, we have
\begin{equation}
e^{-t} \begin{pmatrix} W(t)\\
B(t)
\end{pmatrix} \almostsure
 \lambda_1{\bf v}_1,
\label{Eq:Janson}
\end{equation}
where $\lambda_1$ is the principal eigenvalue of ${\bf A}^\mathsf{T}$ (the eigenvalue of largest real part among all eigenvalues), and ${\bf v}_1$ is the corresponding eigenvector (normalized to length $1$). A straightforward calculation shows 
$$\lambda_1 = 1, \qquad {\bf v}_1 = \begin{pmatrix} 2/3 \\
1/3
\end{pmatrix}.
$$ 
Let $N_i(t)$ be the number of internal nodes (in the original poissonized BST) of degree~$i$, for $i = 1, 2, 3$. A leaf node (of degree 1) corresponds to two white balls in the urn. That is, we have
$$\frac {N_1(t)} {e^t} =  \frac {W(t)} {2e^t}  \almostsure \frac 1 3.$$
An internal node of degree 2 corresponds to an internal node carrying one child (in the original tree) as well as a blue ball in the urn. That is, we have
$$\frac {N_2(t)} {e^t} =  \frac {B(t)} {e^t}  \almostsure \frac 1 3.$$
The size of the tree, $S(t)$, at time $t$ is also random, following a {\em Yule process}~\cite{Simon}, and we have
$$\frac {S(t)}{e^t} \almostsure 1.$$
From this information, we can glean what we need about $N_3(t)$, as we have 
$$ \frac {S(t)}{e^t} = \frac {N_1(t) + N_2(t) + N_3(t)}{e^t} \almostsure 
1.$$
Hence, we have
$$\frac {N_3(t)} {e^t} \almostsure \frac 1 3.$$
We note that, 
for $i=1, 2, 3$, we have $N_i(t) / S(t) \almostsure 1/3$, too, as $t \to \infty$. 
As $N_i(t) / S(t)$ is a proportion, we have $N_i(t) / S(t) \inL 1/3$, as well.
We can now find the relative Gini index of a poissoized BST, $B_{\mathsf P}(t) \in \mathcal{B}_{\mathsf P}(t)$, 
at time $t$ from Equation~(\ref{Eq:binGini}), with the notation adapted
to the situation at hand:
\begin{align*}
G\bigl(B_{\mathsf P}(t)\bigr) &= \frac {N_1(t) N_2(t) + N_2(t) N_3(t) + 2 N_1(t) N_3(t)}
{\E^2 \bigl[S(t)\bigr] \times \E\left[\frac {N_1(t) + 2N_2(t) + 3N_3(t)}
	{S(t)}\right]}
\\ &\almostsure \frac{\frac {N_1(t)} {e^t} \times \frac {N_2(t)} {e^t} 
   + \frac {N_2(t)} {e^t} \times \frac {N_3(t)} {e^t}  + 2 \times\frac {N_1(t)} {e^t} \times \frac {N_3(t)} {e^t} }{ 1/3 + 2/3 + 3/3}
\\ &\almostsure \frac{2}{9},
\end{align*}
as $t \to \infty$.

As proved in~\cite{Domicolo}, the degree Gini index falls in the interval $[0, 1)$. In view of the 
boundedness of $G \bigl(B_{\mathsf P}(t) \bigr)$, the almost-sure convergence implies convergence in the $L_1$ space. That is, 
we have 
$$G^*_{\mathcal{B}_{\mathsf{P}}(t)} = \E \bigl[G \bigl( \mathcal{B}_{\mathsf P}(t) \bigr)\bigr] \to \frac 2 9,$$
as $t\to\infty$. We conclude that the limiting degree-based Gini index of discrete-time BST's and that of the poissonized BST's coincide.  

\section{Binary pyramids}
{\em Binary pyramids}~\cite{Gastwirth1977, Gastwirth1984, Mahmoud1994}, also known as unary-binary trees (a more commonly-used term in computer science and information theory)~\cite{Bacher, Blieberger, Bodini, Flajolet1986}, are models for chain letters, recruitment hierarchy in certain business practices, and similar schemes, with saturation level of 2 for each node. In the context of chain letters, a participant can sell up to a maximum of two letters. The discrete-time pyramid starts out empty at time $n=0$. At time $1$, a node representing the founder of a chain letter scheme appears. 
The founder goes out seeking a purchaser of the letter. When found, the purchaser is linked as a child in the pyramid representing the scheme. Next, these two participants in the scheme compete with equal chance to attract the next purchaser, who is linked as a child in the pyramid to whomever wins the competition. The process proceeds in the following fashion. At each point in discrete time, certain candidate nodes compete with equal chance to attract a new purchaser. As time goes by, the {\em outdegree}\footnote{The outdegree of a node is the number of edges emanating out from it.} of a node can increase from 0 to 1 or 2. Once the outdegree of a node reaches the cap~2, the node is considered saturated, having reached the prespecified quota of 2. In other words, such a node is taken out of the competition to acquire new purchasers in subsequent rounds.

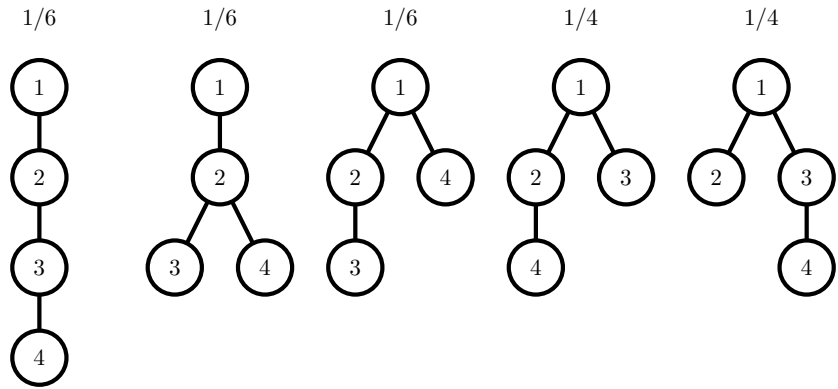
\begin{figure}[thb]
	\begin{center}
		\begin{tikzpicture}[scale=0.6]
		\coordinate (A) at (-5,7);
		\coordinate (B) at (-5,5);
		\coordinate (C) at (-5,3);
		\coordinate (D) at (-5,1);

		\draw [ultra thick] (A)--(B)--(C)--(D);
		\draw [ultra thick,fill=white] (-5,7) circle [radius=0.6];
		\draw [ultra thick,fill=white] (-5,5) circle [radius=0.6];
		\draw [ultra thick,fill=white] (-5,3) circle [radius=0.6];
		\draw [ultra thick,fill=white] (-5,1) circle [radius=0.6];
		
		\node[draw=white, scale = 0.8] at (-5, 7) {$1$};
		\node[draw=white, scale = 0.8] at (-5, 5) {$2$};
		\node[draw=white, scale = 0.8] at (-5, 3) {$3$};
		\node[draw=white, scale = 0.8] at (-5, 1) {$4$};
		\node[draw=white, scale = 0.8] at (-5, 8.5) {$1/6$};
		
		\coordinate (A) at (-1,7);
		\coordinate (B) at (-1,5);
		\coordinate (C) at (-2,3);
		\coordinate (D) at (0,3);

		\draw [ultra thick] (A)--(B);
		\draw [ultra thick] (B)--(C);
		\draw [ultra thick] (B)--(D);
		\draw [ultra thick,fill=white] (-1,7) circle [radius=0.6];
		\draw [ultra thick,fill=white] (-1,5) circle [radius=0.6];
		\draw [ultra thick,fill=white] (-2,3) circle [radius=0.6];
		\draw [ultra thick,fill=white] (0,3) circle [radius=0.6];
		
		\node[draw=white, scale = 0.8] at (-1, 7) {$1$};
		\node[draw=white, scale = 0.8] at (-1, 5) {$2$};
		\node[draw=white, scale = 0.8] at (-2, 3) {$3$};
		\node[draw=white, scale = 0.8] at (0, 3) {$4$};
		\node[draw=white, scale = 0.8] at (-1, 8.5) {$1/6$};

		\coordinate (A) at (3,7);
		\coordinate (B) at (2,5);
		\coordinate (C) at (2,3);
		\coordinate (D) at (4,5);
		
		\draw [ultra thick] (A)--(B)--(C);
		\draw [ultra thick] (A)--(D);
		\draw [ultra thick,fill=white] (3,7) circle [radius=0.6];
		\draw [ultra thick,fill=white] (2,5) circle [radius=0.6];
		\draw [ultra thick,fill=white] (2,3) circle [radius=0.6];
		\draw [ultra thick,fill=white] (4,5) circle [radius=0.6];
		
		\node[draw=white, scale = 0.8] at (3, 7) {$1$};
		\node[draw=white, scale = 0.8] at (2, 5) {$2$};
		\node[draw=white, scale = 0.8] at (2, 3) {$3$};
		\node[draw=white, scale = 0.8] at (4, 5) {$4$};
		\node[draw=white, scale = 0.8] at (3, 8.5) {$1/6$};

		\coordinate (A) at (7,7);
		\coordinate (B) at (6,5);
		\coordinate (C) at (8,5);
		\coordinate (D) at (6,3);
		
		\draw [ultra thick] (A)--(B)--(D);
		\draw [ultra thick] (A)--(C);
		\draw [ultra thick,fill=white] (7,7) circle [radius=0.6];
		\draw [ultra thick,fill=white] (6,5) circle [radius=0.6];
		\draw [ultra thick,fill=white] (8,5) circle [radius=0.6];
		\draw [ultra thick,fill=white] (6,3) circle [radius=0.6];
		
		\node[draw=white, scale = 0.8] at (7, 7) {$1$};
		\node[draw=white, scale = 0.8] at (6, 5) {$2$};
		\node[draw=white, scale = 0.8] at (8, 5) {$3$};
		\node[draw=white, scale = 0.8] at (6, 3) {$4$};
		\node[draw=white, scale = 0.8] at (7, 8.5) {$1/4$};
		
		\coordinate (A) at (11,7);
		\coordinate (B) at (10,5);
		\coordinate (C) at (12,5);
		\coordinate (D) at (12,3);
		
		\draw [ultra thick] (A)--(C)--(D);
		\draw [ultra thick] (A)--(B);
		\draw [ultra thick,fill=white] (11,7) circle [radius=0.6];
		\draw [ultra thick,fill=white] (10,5) circle [radius=0.6];
		\draw [ultra thick,fill=white] (12,5) circle [radius=0.6];
		\draw [ultra thick,fill=white] (12,3) circle [radius=0.6];
		
		\node[draw=white, scale = 0.8] at (11, 7) {$1$};
		\node[draw=white, scale = 0.8] at (10, 5) {$2$};
		\node[draw=white, scale = 0.8] at (12, 5) {$3$};
		\node[draw=white, scale = 0.8] at (12, 3) {$4$};
		\node[draw=white, scale = 0.8] at (11, 8.5) {$1/4$};
		
		\end{tikzpicture}
	\end{center}
	\caption{Binary pyramids of order 4 with their probabilities at the top.}
	\label{Fig:pyramids}
\end{figure}

Basically, a binary pyramid is a rooted binary tree, on which another nonuniform distribution imposed. 
The five possible binary pyramids emerging at time $4$ (after four insertions) are shown in Figure~\ref{Fig:pyramids}, together with their probabilities 
above each. The four nodes in the leftmost pyramid in Figure~\ref{Fig:pyramids} are all unsaturated and can recruit the next purchaser---there is an insertion position under each node (a total of four insertion positions). In the second pyramid from the left in Figure~\ref{Fig:pyramids}, the node labeled with~$2$ is saturated as it has reached the cap; the other three nodes are unsaturated---there is an insertion position under each of the nodes respectively labeled with 1, 3 and~4 (a total of three insertion positions). Based on this observation, we notice that binary pyramids of the same order may not have the same number of insertion positions which would appear as external nodes in the extended graph.

The randomness of the insertion positions is in contrast with BST's determinism, where all trees of the same order have the same number of external nodes. This has significantly hampered research on pyramids in the past. They are not as well investigated as BST's or random recursive trees~\cite{Mahmoud1991, Mahmoud1992}, another class of random graph models of substantial interest.   
 
\subsection{Gini index of discrete-time binary pyramids} 
We can establish correspondence between binary pyramids and two-color \Polya\ urns. Toward this end, we introduce another color code for the external nodes (insertion positions) in extended binary pyramids. For unsaturated nodes in the original binary pyramid, external 
nodes connected to candidate nodes of outdegree 0
(leaf nodes in the original tree) 
are painted white, while those connected to candidate nodes of outdegree 1 (internal nodes in the original tree) are painted blue. Next, think of the colored external nodes as balls in an urn.

\begin{figure}[tbh]
	\begin{center}
		\begin{minipage}{0.48\textwidth}
			\begin{tikzpicture}[scale=0.67]
			\coordinate (A) at (-3,4);
			\coordinate (B) at (3,4);
			
			\coordinate (D) at (-3,2);
			\draw [ultra thick] (A)--(D);
			\draw [ultra thick,fill=white] (-3,4) circle [radius=0.6];
			\node at (-3,2) [rectangle, draw, ultra thick, scale = 2.5, fill=white](){};

			\draw [->,line width=1pt] (-1,2) -- (0,2); 
			
			\coordinate (F) at (2,2);
			\coordinate (G) at (4,2);
			\coordinate (H) at (2,0);
			\draw [ultra thick] (B)--(F);
			\draw [ultra thick] (B)--(G);
			\draw [ultra thick] (F)--(H);
			\draw [ultra thick,fill=white] (3,4) circle [radius=0.6];
			\draw [ultra thick,fill=white] (2,2) circle [radius=0.6];
			\node at (4,2) [rectangle, draw, ultra thick, scale = 2.5, fill=blue](){};
			\node at (2,0) [rectangle, draw, ultra thick, scale = 2.5, fill=white](){};
			\end{tikzpicture}
		\end{minipage}
		\hfill
		\begin{minipage}{0.48\textwidth}
			\begin{tikzpicture}[scale=0.67]
			\coordinate (A) at (-3,4);
			\coordinate (B) at (3,4);
			
			\coordinate (D) at (-4,2);
			\coordinate (E) at (-2,2);
			\draw [ultra thick] (A)--(D);
			\draw [ultra thick] (A)--(E);
			\draw [ultra thick,fill=white] (-3,4) circle [radius=0.6];
			\draw [ultra thick,fill=white] (-4,2) circle [radius=0.6];
			\node at (-2,2)[rectangle, draw, ultra thick, scale = 2.5, fill=blue](){};
			
			\draw [->,line width=1pt] (-1,2) -- (0,2); 
			
			\coordinate (F) at (2,2);
			\coordinate (G) at (4,2);
			\coordinate (H) at (4,0);
			\draw [ultra thick] (B)--(F);
			\draw [ultra thick] (B)--(G);
			\draw [ultra thick] (G)--(H);
			\draw [ultra thick,fill=white] (3,4) circle [radius=0.6];
			\draw [ultra thick,fill=white] (2,2) circle [radius=0.6];
			\draw [ultra thick,fill=white] (4,2) circle [radius=0.6];
			\node at (4,0) [rectangle, draw, ultra thick, scale = 2.5, fill=white](){};
			\end{tikzpicture}
		\end{minipage}
		\caption{
		Correspondence of a two-color P\'{o}lya urn to the evolution of  an extended binary pyramid.}
		\label{Fig:polyaurnpyramid}
	\end{center}
\end{figure}
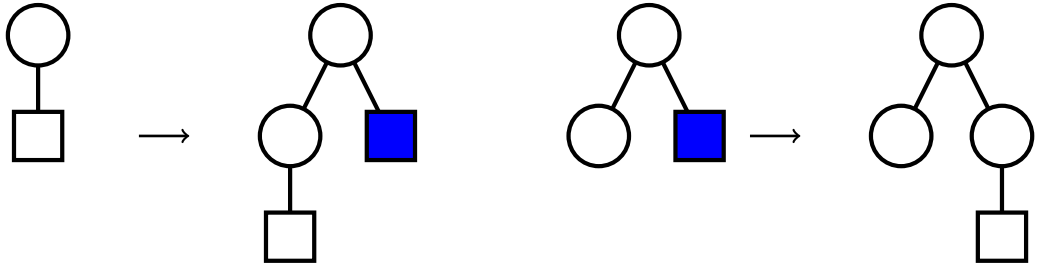

When a white external node (white ball in the urn) is picked, an insertion position under a leaf node is converted into a new leaf node with a white external node under it. An additional blue external node appears as a sibling of the new leaf node. The net gain is one blue ball. Alternatively, if we pick a blue external node, the corresponding insertion position is converted into an internal node (leaf node), with one white external node underneath it. The parent of the lost blue external node becomes saturated. Note that there are no external nodes under saturated nodes, and they do not correspond to any balls in the urn. See Figure~\ref{Fig:polyaurnpyramid} for a graphic interpretation. 

The ball addition mechanism can be represented by the following replacement matrix:
$${\bf A} = \begin{pmatrix} 0  & 1\\
1   & -1
\end{pmatrix}, $$
where white and blue are indexed from top to bottom, and from left to right. Let $W_n$ and $B_n$ be the number of white and blue balls in the urn at time~$n$, respectively. The ball replacement matrix has the principal eigenvalue $\lambda_1 =  (\sqrt 5 - 1)/2$, and a corresponding ($L_1$-normalized) principal eigenvector
${\bf v_1} =\begin{pmatrix} (\sqrt 5 - 1)/2\\
(3-\sqrt 5)/2 \end{pmatrix}$.

Again, appealing to the urn theory in~\cite{Janson}, we have
$$\frac 1 n \begin{pmatrix} W_n\\
B_n
\end{pmatrix} 
\almostsure \lambda_1{\bf v}_1
= \begin{pmatrix}
(3 - \sqrt{5})/2
\\ \sqrt{5} - 2
\end{pmatrix}.$$

Let $N_i(n)$ be the number of internal nodes 
of degree $i$ in the original binary pyramid (not yet extended), for $i = 1, 2, 3$. A leaf node (of degree 1 and of outdegree 0) corresponds to one white ball in the urn. We thus have
$$\frac {N_1(n)} n =  \frac {W_n} {n}  \almostsure \frac {3 - \sqrt5}{2}.$$
An internal node of degree 2 (and of outdegree 1) corresponds to a blue ball in the urn. Therefore,
we obtain
$$\frac {N_2(n)} n =  \frac {B_n} n \almostsure   \sqrt 5-2.$$ 
After $n$ insertions in an empty pyramid, 
there are $n$ nodes. Thus, the number of saturated nodes is $N_3(n) = n - N_1(n) - N_2(n)$, and we have
$$\frac{N_3(n)}{n} \almostsure 1  - \frac{3 - \sqrt 5}{2} -   (\sqrt 5 -2)
= \frac{3 - \sqrt{5}}{2}.$$ Since $N_i(n)/n$ is a proportion, for $i=1,2,3$,
these almost-sure relations hold in $L_1$, too.
Hence, the degree Gini index of a binary pyramid, $P(n)$, of order $n$, is 
\begin{align*}
G\bigl(P(n)\bigr) &= \frac {N_1(n) N_2(n) + N_2(n) N_3(n) + 2 N_1(n) N_3(n)}
{\displaystyle n^2 \,\E\left[\frac {N_1(n) + 2N_2(n) + 3N_3(n)}
	n\right]} \\
&=\frac {\frac {N_1(n)} n\times \frac{N_2(n)} n + \frac {N_2(n)} n
   \times \frac {N_3(n)} n + 2\,  \frac {N_1(n)} n \times \frac {N_3(n)} n}
{ \left(\E \bigl[N_1(n)\bigr] + 2\E\bigl[N_2(n)\bigr] + 3\E\bigl[N_3(n)\bigr]\right)/n}
\end{align*}
By the claimed $L_1$ convergence, we further have      
\begin{align*}
G(P_n) &\almostsure \left(\frac{(3 - \sqrt{5})(\sqrt{5} - 2)}{2} + \frac{(\sqrt{5} - 2)(3 - \sqrt{5})}{2} + \frac{(3 - \sqrt{5})^2}{2}\right)\\
&\qquad \qquad{} \times \left(\frac{3 - \sqrt{5}}{2} + 2(\sqrt{5} - 2) + \frac{3(3 - \sqrt{5})}{2}\right)^{-1} \\
&= \sqrt 5 - 2\\
&\approx 0.236068.
\end{align*}
In view of the 
boundedness of $G(P_n)$, the almost-sure convergence of the Gini index
of a pyramid implies convergence of the mean, too. That is, 
for the class of binary pyramids of size $n$, denoted by ${\cal P}(n)$, we have 
$$G^*_{{\cal P}(n)} =\E \bigl[G(P_n)\bigr] \to \sqrt 5 - 2 \approx 0.236068,$$
as $n \to \infty$.
\subsection{Gini index of poissonized binary pyramids}
The poissonized binary pyramid has 
the same underlying \Polya\ urn, with the renewals occurring at an ever-accelerating rate (in view of an increasing number of $\Exp(1)$
interarrival times). The computations are essentially the same, mutatis mutandis, as in BST's (except for a different replacement matrix). We present only the result, omitting details of the derivation. For the class of poissonized binary pyramids at time $t$, denoted ${\cal P_{\mathsf{P}}}(t)$, the degree-based Gini index is
$$G^*_{{\cal P_{\mathsf{P}}} (t)} \to \sqrt 5 - 2 \approx 0.236068,$$
as $t \to \infty$, exactly the same as that of the discrete-time binary pyramids.
\section{Uniform caterpillars} 
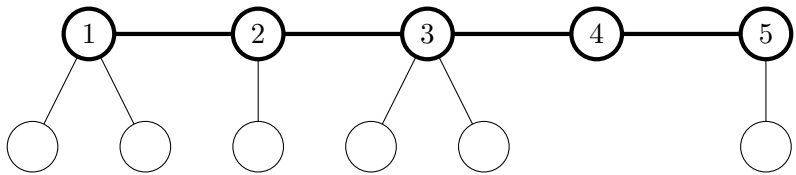
\begin{figure}[tbh]
	\begin{center}
		\begin{tikzpicture}[scale=3]
		\coordinate (A) at (-1.25, 0);
		\coordinate (B) at (-0.50, 0);
		\coordinate (C) at (0.25, 0);
		\coordinate (D) at (1.25, 0);
		\coordinate (E) at (1.75, 0);
		\draw [ultra thick] (A)--(E);
		
		\coordinate (F) at (-1.5, -0.5);
		\coordinate (G) at (-1.0, -0.5);
		\coordinate (H) at (-0.5, -0.5);
		\coordinate (I) at (0.0, -0.5);
		\coordinate (J) at (0.5, -0.5);
		\coordinate (K) at (1.75, -0.5);
		\draw  (A)--(F);
		\draw  (A)--(G);
		\draw  (B)--(H);
		\draw  (C)--(I);
		\draw  (C)--(J);
		\draw  (E)--(K);
		\draw
		(-1.5, -0.5) node [circle=0.1,draw,fill=white] {\phantom{$1$}}
		(-1, -0.5) node [circle=0.1,draw,fill=white] {\phantom{$2$}}
		(-0.5, -0.5) node [circle=0.1,draw,fill=white] {\phantom{$3$}}
		(0, -0.5) node [circle=0.1,draw,fill=white] {\phantom{$4$}}
		(0.5, -0.5) node [circle=0.1,draw,fill=white]{\phantom{$5$}}
		(1.75, -0.5) node [circle=0.1,draw,fill=white] {\phantom{$6$}};
		\draw
		(-1.25, 0) node [ultra thick,circle=0.1,draw,fill=white] {{$1$}}
		(-0.5, 0) node [ultra thick,circle=0.1,draw,fill=white] {{$2$}}
		(0.25, 0) node [ultra thick,circle=0.1,draw,fill=white] {{$3$}}
		(1, 0) node [ultra thick,circle=0.1,draw,fill=white] {{$4$}}
		(1.75, 0) node [ultra thick,circle=0.1,draw,fill=white] {{$5$}};
		\end{tikzpicture}
		\caption{A caterpillar of a central path of five nodes (thickened boundaries)
			and six attachments.}
		\label{Fig:caterpillar}
	\end{center} 
\end{figure}

The reference~\cite{Harary} poetically describes a caterpillar as ``a tree which metamorphoses into a path, when its cocoon of endpoints is removed.'' In graph theory, a {\em caterpillar} is a tree such that every node has distance at most $1$ from a {\em central path}, which is also called a {\em spine} in the literature. An example of a caterpillar is shown in Figure~\ref{Fig:caterpillar}.
There are different names for the vertices at distance 1 from the spine, such as monovalent vertices and endpoints. For simplicity, we call them {\em attachments} in the rest of this manuscript.

In chemistry literature, caterpillars are called Gutman trees or benzenoid trees; see~\cite{ElBasil1990}.
A topological perspective of caterpillars is given in~\cite{Gutman, Schmuck}. Caterpillars can model the growth of algae on a food source. The model is also pertinent to the formation of computer networks with one common ``bus'' or ``hub'' of servers, to which users subscribe to receive Internet services.

Random caterpillars can evolve in different ways. Our goal is to probe the duality arising from embedding in continuous time. Therefore, we choose a growth model that can be poissonized. 

We consider {\em uniform} caterpillars: In the discrete-time version, initially (at time $n=0$) we have a spine consisting of a fixed number of nodes, say $s$. At every subsequent point $n\ge 1$ of discrete time, a node is attached to the spine by choosing a spine node uniformly at random (all spine nodes have equal probability, which is namely $1/s$). A distance-based Gini index for this class of random caterpillars is studied in~\cite{Balaji}.

A poissonized version of the discrete-time uniform counterpart grows as follows. Each spine node has a renewal process associated with it. The renewals occur at independent interarrival times that are distributed like $\Exp(1)$ random variable. When a renewal occurs at a spine node, a new node is attached to that spine node, and the stochastic process of renewals continues in a memoryless fashion, in view of the exponential interarrival times.  According to the setup, each spine node generates a {\em point Poisson process}~\cite{Ross}; by time $t$, the number of nodes attached to a spine node has a Poisson distribution with parameter $t$. 

We shall call the class of caterpillars that grow in discrete time ${\cal C}(n)$, and the poissonized counterpart ${\cal C_\mathsf{P}}(t)$. We investigate the degree-based Gini index proposed in~\cite{Domicolo} for both ${\cal C}(n)$ and ${\cal C_\mathsf{P}}(t)$.

\subsection{Gini index of discrete-time uniform caterpillars}
\label{Sec:discaterpillar} 
Suppose that we are given a caterpillar with a spine which consists of a finite number of $s\ge 1$ nodes. 
Let us label the $s$ spine nodes from one end to the other sequentially with distinct numbers from the set $\{1, 2, \ldots, s\}$. Let $X_i(n)$ be the number of nodes attached to the spine node labeled with $i$ by discrete time $n$. At that time, the size of the caterpillar is fixed; it is $\sum_{i = 1}^{s} X_i(n) + s = n + s$. The expected degree of a randomly selected node in a randomly generated caterpillar, $C(n) \in \mathcal{C}(n)$, is
\begin{align*}
\E\left[D^* \bigl(C(n)\bigr)\right] &= \frac{\bigl(X_{1}(n) + 1\bigr) + \sum_{i = 2}^{s - 1} \bigl(X_{i}(n) + 2\bigr) + \bigl(X_{s}(n) + 1\bigr) + \sum_{i = 1}^{s}X_{i}(n)}{n + s}
\\&= \frac{2n +2s-2} {n+s}\\
&= 2 + O(1/n);
\end{align*}
in the numerator of this formula, the first three terms account for the contributions of spine nodes (with a special handling of its two end points), and the last term accounts for the contributions of the attachments.

To execute Formula~(\ref{Eq:singlegraph}), we need to compute the absolute differences in the numerator. We split the overall contributions to the absolute differences into inter-spine contributions, spine-attachments contributions, and attachment-attachment contributions, where the latter is 0, as all attachments have degree $1$. 

\subsubsection{Inter-spine contributions}
With the exception of the two end vertices on the spine, each other spine node has degree equal to the number of nodes attached to it plus $2$. The degree of a node at one of the ends of the spine is equal to the number of nodes attached to it plus $1$. Hence, the contributions of inter-spine nodes to the numerator of the degree-based Gini index in Equation~(\ref{Eq:singlegraph})
is the expected value of 
\begin{align*}
&\sum_{1< i<j<s} \big|X_i(n) - X_j(n)\big|  +  \sum_{i=2}^{s-1} \big|X_i(n) + 1 - X_1(n)\big|  +  \sum_{i=2}^{s-1} \big|X_i(n) + 1 - X_s(n)\big| \\
&\qquad\qquad {} +  \big |X_1(n) - X_s(n)\big|.
\end{align*}
In view of the identical distribution of the attachments, this expectation is 
\begin{align*}
&{s - 2 \choose 2} \E\big|X_2(n) - X_1(n)\big|  +  2(s - 2) \E \big|X_2(n) + 1 - X_1(n)\big| + 
\\ &\qquad\qquad{}+ \E\big |X_1(n) - X_s(n)\big|.
\end{align*}
Each of these terms is $O(n)$, since $s$ is assumed finite. Hence, the expectation of entire inter-spine contributions is also $O(n)$.
\subsubsection{Spine-attachment contributions}
We compute the spine-attachment contributions by handling the two nodes at the ends of the spine separately. The expectation of spine-attachment contributions is given by
\begingroup
\allowdisplaybreaks
\begin{align}
&\E\Bigl[n \big| \bigl(X_1(n) + 1 \bigr)  - 1 \big| +  n \sum_{i=2}^{s-1} \big| \bigl(X_i(n) +2 \bigr)-1 \big| + n \, \big| \bigl(X_s(n) + 1 \bigr) - 1\big |\Bigr] \nonumber\\
&\qquad{} = n \, \E\bigl[X_1(n)\bigr]  + n \sum_{i=2}^{s-1} \bigl(\E\bigl[X_i(n)\bigr] +1\bigr)  +  n \, \E\bigl[X_s(n)\bigr]\nonumber
\\&\qquad{} =  n \sum_{i=1}^s X_i(n)  + (s-2)n \label{Eq:numerator}\\
&\qquad = n^2 + (s - 2)n. \nonumber
\end{align}
\endgroup

Putting these elements together, we find the degree-based Gini index of the class of random uniform caterpillars after $n$ insertions to be
$$G^*_{{\cal C}(n)} = \frac {n^2 + O(n)}  {(n+s)^2 (2 + O(1/n))} \to \frac 1 2,$$
as $n \to \infty$.
\subsection{Gini index of poissonized uniform caterpillars}
\label{Sec:poiscaterpillar}
We denote a Poisson random variable with parameter $\lambda$ by
Poi$(\lambda)$. 
As in the discrete-time uniform caterpillar, let us label the $s$ spine nodes of a poissonized uniform caterpillar from one end to the other sequentially with the distinct numbers in set $\{1, 2, \ldots, s\}$.
Let $Y_i(t)$ be the number of nodes attached to the spine node labeled with $i$ by time~$t$. Recall that the number of attachments to a spine node is a Poisson process with intensity $t$. 
Thus, the size of a poissonized uniform caterpillar by time $t$, denoted by~$S(t)$, is random. Namely, it is
$$S(t) = s + \sum_{i=1}^s Y_i(t),$$
where the sum is that of $s$ independent {\rm Poi($t$)} random variables. So, it is itself a Poisson random variable with parameter $st$. The average size by time $t$ is
\begin{equation}
\E\bigl[S(t)\bigr]= s + st.
\label{Eq:size}
\end{equation} 
We can now compute the average degree of a randomly chosen node in a poissonized caterpillar, $C_{\mathsf{P}}(t)\in \mathcal{C}_{\mathsf{P}}(t)$, by time $t$:
\begin{align*}
\E\bigl[D^* \bigl(C_{\mathsf{P}}(t)\bigr)\, | \,  S(t)\bigr] &= \frac{\bigl(Y_{1}(t) + 1\bigr) + \sum_{i = 2}^{s - 1} \bigl(Y_{i}(t) + 2\bigr) + \bigl(Y_{s}(t) + 1\bigr) + \sum_{i = 1}^{s}Y_{i}(t)}{S(t)}
\\
&= \frac{2S(t) -2} {S(t)}\\
&= 2 - \frac{ 2}{S(t)}.
\end{align*}
Recall that $\bigl(S(t) - s\bigr)$ has a Poisson distribution with parameter $st$. Given any sufficiently large finite number $M > 0$, we have
$$\lim_{t \to \infty} \Prob\bigl(S(t) - s < M\bigr) = \lim_{t \to \infty} e^{-st} \sum_{k = 0}^{\lfloor M \rfloor} \frac{(st)^k}{k!} = 0,$$
and we conclude that $2/ S(t)$ converges to $0$ in probability. Besides, with $2/S(t)$ uniformly bounded (by
$2$), this in-probability convergence leads to
convergence on average. Hence, we get
$$\E\left[D^*\bigl(C_{\mathsf{P}}(t)\bigr)\right] = 2 + o(1),$$
as $t \to \infty$. 

The inter-spine and spine-attachments relations can be computed in the same way as in the discrete-time uniform caterpillar, so we get an analogue of Equation~(\ref{Eq:numerator}), with $X_i(n)$ replaced by $Y_i(t)$ and the multiplier $n$ replaced by the expected number of
attachments, 
given in~(\ref{Eq:size}). The numerator of the Gini index becomes
$$ \E\bigl[\Poi^2(st) + (s-2)\Poi(st) + O(t)\bigr] = (st + s^2 t^2) + (s - 2)st + O(t) = s^2 t^2 + O(t). $$ 

We now have all the elements for the calculation of the degree-based Gini index for the class of poissonized uniform caterpillars, $\mathcal{C}_{\mathsf P}(t)$:
$$G^*_{\mathcal{C}_{\mathsf P}(t)} = \frac {s^2 t^2 + O(t)} {(s+st)^2(2 + o(1))} 
= \frac 1 2 + O\left(\frac{1}{t}\right) \to \frac{1}{2},$$
as $t \to \infty$.  
This is indeed identical to the limiting degree-based Gini index of the discrete uniform caterpillars.

\subsection{Another version of Gini index}
A slightly different version of Gini index is proposed in~\cite{Zhang}, where all spine nodes form a population, and the number of attachments linked to each spine node is thought of as its wealth. 
This version of Gini index of a uniform caterpillar at time $n$ is concerned only with the inter-spine absolute differences, as a measure of inequality of wealth among the spine nodes. Let the {\em wealth} 
(the number of attachments) of the $i$th node on the spine be $W_i$.
The authors of~\cite{Zhang} define an inequality measure by the following formula
\begin{equation}
\tilde{G}^{*}_{\mathcal{C}(n)} = \frac{\sum_{i = 1}^{s} \sum_{j = 1}^{s} \E \bigl|W_i(n) - W_j(n)\bigr|}{2 s \sum_{i = 1}^{s} W_i(n)} = \frac{(s^2 - s)\, \E\bigl|W_1(n) - W_2(n)\bigr|}{2 s n}, \label{Eq:anotherGini}
\end{equation}
where $W_i(n)$'s follow a multinomial distribution with the number of trails $n$ 
and event probabilities $\{1/s, 1/s, \ldots, 1/s\}$. 
Thus, each of $W_1(n)$ and $W_2(n)$ has a binomial distribution with parameters $n$ and $1/s$. We have $W_1/n \to 1/s$
and $W_2/n \to 1/s$ in the $L_1$ space and $(W_1 - W_2)/n \to 0$ in $L_1$ as well. According to the {\em Continuous Mapping Theorem}, 
we get $|W_1 - W_2|/n \to 0$ in the $L_1$ space, too. Hence, we conclude that, as long as $s = o(n)$, the Gini index $\tilde{G}^{*}_{\mathcal{C}(n)} \to 0$ when $n \to \infty$. 

This version of Gini index also converges to $0$ for poissonized uniform caterpillars. The proof can be done by establishing analogous arguments, while the terms $W_i(n)$'s in the numerator of Equation~(\ref{Eq:anotherGini}) are replaced by $\Poi(1)$ random variables, say $Q_i(t)$'s, and the term $\sum_{i = 1}^{s} W_i(n)$ in the denominator is replaced by $\sum_{i = 1}^{s}\E \bigl[Q_i(t)\bigr] = st$. Thus, we conclude that the duality principle carries over to this particular version of Gini index for discrete-time caterpillars and the poissonized counterpart.
\section{Preferential attachment caterpillars}
In~\cite{Szym}, a novel mechanism of edge emergence is proposed for growth in network. This method is called {\em preferential attachment}. The primary feature of this attachment scheme is that a node of larger degree has a higher probability to attract newcomers as time goes by. 
The model was very widely popularized after~\cite{Barabasi} connected it
to scale-free networks of real-world networks following power laws.
 
Statisticians exploit preferential attachment to mathematically interpret phenomena from various disciplines, such as the {\em Matthew Effect} in economics~\cite{Matthew}, the feature of citation networks in social science~\cite{desolla}, and the expansion of the Internet in computer science~\cite{Barabasi}. 

A class of preferential attachment caterpillars is proposed in~\cite{Zhang}. 
We shall call this class $\mathcal{C}^{\mathsf{PA}}(n)$
At time 0, start with a spine consisting of $s$ nodes. At each timestamp $n \ge 1$, the probability of selecting a spine node for the next attachment
is proportional to its degree at time $n - 1$. We reuse the notation, $X_i(n)$, 
to denote the number of attachments linked to the spine node labeled with $i$ (from left to right) at time $n$. The probability that spine node $i$ is chosen to attract a newcomer at time $(n + 1)$ is $(X_i(n) + 1)/\left(\sum_{i = 1}^{s} X_i(n) + 2s - 2\right)$  for $i = 1$ and $i = s$; the probability is $(X_i(n) + 2)/\left(\sum_{i = 1}^{s} X_i(n) + 2s - 2\right)$  for $i = 2, \ldots, s - 1$. The expected degree of a randomly selected node in a random caterpillar in this class, $C^{\mathsf{PA}}(n) \in \mathcal{C}^{\mathsf{PA}}(n)$, is identical to that for the uniform caterpillar class; that is,
$$\E\left[D^* \bigl(C^{\mathsf{PA}}(n)\bigr)\right] =\frac{2n +2s-2} {n+s}.$$

To compute the numerator of Equation~(\ref{Eq:singlegraph}) for this class, we reconsider three kinds of contributions: inter-spine contributions, spine-attachment contributions, and attachment-attachment contributions, where the last add $0$.
Establishing analogous arguments, we find that the expectations of the inter-spine contributions and the spine-attachment contributions are respectively $O(n)$ and $n^2 + (s - 2)n$. We 
conclude that
$$G^*_{\mathcal{C}^{\mathsf{PA}}_n} = \frac{n^2 + O(n)}{n^2 \times (2n + 2s - 2)/(n + s)}\to \frac{1}{2},$$
as $n \to \infty$. The limiting degree-based Gini index of poissonized preferential attachment caterpillars is also $1/2$, where the computation is identical to that for the uniform class. So, we do not repeat ourselves in the article. The primary reason of the identicality
of computation is that the degree-based Gini index for random caterpillars is based on the totality of attachments, but invariant with respect to  the evolutionary characteristics. In conclusion, the duality theory is verified for preferential attachment caterpillars growing in discrete-time and the poissonized counterparts.

\section{Conclusion}
The duality theory, in its general sense, argues that a problem can be approached from two perspectives, and these share certain properties. In providing the previous examples and proving the connection between the degree-Gini index of discrete-time and poissonized trees, we are looking to move closer to a strong understanding of the duality theory. Our examples for certain broad classes of trees (binary search trees, binary pyramids and random caterpillars) provide proof that the degree-Gini index of a class of discrete-time trees is asymptotically the same as the degree-Gini index for the same class of trees that have undergone poissonization.

Having this initial proof of the strength of the duality theory in relation to the degree-Gini index, there are future paths that can be taken to continue work on this topic. A direct extension to be pursued is the calculation and comparison of the degree-Gini indices of additional classes of graphs. This information would help to improve the understanding and acceptance of the duality theory for this type of problem.

\end{document}